\documentclass[12pt,reqno]{article}
\usepackage[usenames]{color}
\usepackage[colorlinks=true,
linkcolor=webgreen,
filecolor=webbrown,
citecolor=webgreen]{hyperref}
\definecolor{webgreen}{rgb}{0,.5,0}
\definecolor{webbrown}{rgb}{.6,0,0}

\usepackage{epsfig}                                                                               
\usepackage{amsmath,amssymb,amsthm,eucal}
\usepackage{amsfonts}

\newtheorem{theorem}{Theorem}

\newenvironment{packed_enumerate}{
\setlength{\parsep}{0pt}
\setlength{\parskip}{0pt}
\begin{enumerate}
  \setlength{\itemsep}{1pt}
  \setlength{\parsep}{0pt}
  \setlength{\parskip}{0pt}
}{\end{enumerate}}

\begin{document}

\begin{center}
\vskip 1cm{\LARGE\bf
Designing peg solitaire puzzles}
\vskip 1cm
{\large George I. Bell}\\
Boulder, CO 80303, USA\\
\href{mailto:gibell@comcast.net}{\tt gibell@comcast.net} \\
\end{center}

\vskip .2 in
\begin{abstract}
Peg solitaire is an old puzzle with a 300 year history.
We consider two ways a computer can be utilized to find interesting peg solitaire puzzles.
It is common for a peg solitaire puzzle to begin from a symmetric board position,
we have computed solvable symmetric board positions for four board shapes.
A new idea is to search for board positions which have a unique starting jump leading to a solution.
We show many challenging puzzles uncovered by this search technique.
Clever solvers can take advantage of the uniqueness property to help solve these puzzles.
\end{abstract}

\section{Introduction}
\label{sec:intro}

Peg solitaire was invented in France in the late 17th century,
where it started an early puzzle craze.
Today most people recognize the puzzle, although its popularity has declined.

\begin{figure}[htb]
\centering
\epsfig{file=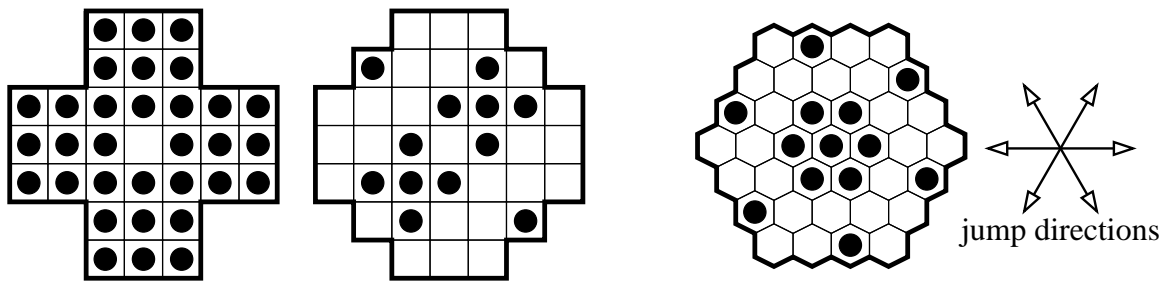}
\caption{Sample puzzles on the 33-hole ``English" board, the 37-hole ``French" board and the 37-hole hexagon board.}
\label{fig1}
\end{figure}

We will refer to a board location as a {\bf hole}, which can either be empty
or occupied by a {\bf peg}.
Figure~\ref{fig1} shows three peg solitaire boards---the
first two boards are based on a square lattice of holes,
while the third is based on a triangular lattice.
While the first two boards are common,
the 37-hole hexagon board is not.
Pressman Toy Company has manufactured this board under the name {\it Think 'N Jump},
although it is not identical since they removed some of the outer jumps
(one can still play on this board by allowing these jumps).

The puzzle begins from some specified pattern of pegs,
three examples are shown in Figure~\ref{fig1}.
A {\bf jump} consists of one peg jumping a neighbor into an empty hole, the jumped peg is removed.
Jumps are allowed only along lattice lines, i.e. along columns and rows for the first two boards,
and along three lines on the hexagon board (as indicated by the arrows).
The goal of the puzzle is to choose a sequence of jumps which finish at a board position with one peg.
The hole where the final peg ends at is called a {\bf finishing hole}.
A board position where there exists a sequence of jumps ending with one peg
is said to be {\bf solvable}.

Figure~\ref{fig1}a shows the starting board position of a special puzzle called the {\bf central game}.
The starting position has square symmetry, and the goal is to finish with one peg in the center,
a board position which not only has square symmetry, but is also the {\bf complement} of
the starting position (where each peg is replaced by a hole, and vice-versa).

If we take either of the 37-hole boards, and fill them with pegs, but remove the central peg,
then we are in an unsolvable board position (we will prove this).
Therefore, the analogous ``central game" is unsolvable on these two boards,
but many other symmetric board states are solvable.
The configuration shown in Figure~\ref{fig1}b, for example, is solvable.
This board position is symmetric with respect to reflection about both diagonals.
Finally, the board position in Figure~\ref{fig1}c is solvable
and is symmetric with respect to 60$^{\circ}$ rotations.
If you solve this puzzle you will discover that the finishing hole
is not the central hole.

For the English and French boards, there are a total of seven symmetries possible for a
configuration of pegs, summarized in Table~\ref{tab1}.
These correspond to subgroups of $D_8$, the dihedral group of
order 8 (the symmetries of the square).
If the mirror is parallel to lattice lines, we call it an
``orthogonal reflection", otherwise it is a ``diagonal reflection".
This distinction is obvious on the square lattice boards,
but more subtle on the 37-hole hexagon board.

\begin{table}[htb]
\begin{center} 
\begin{tabular}{ c  l  c  c }
Type & Symmetry description & Order & Examples \\
\hline
1 & square symmetry & 8 & Fig. \ref{fig1}a, \ref{fig12}a  \\
2 & 90$^{\circ}$ rotation & 4 & Fig. \ref{fig5}c, \ref{fig12}b \\
3 & both diagonal reflections & 4 & Fig. \ref{fig1}b \\
4 & both orthogonal reflections & 4 & Fig. \ref{fig5}a, \ref{fig12}c \\
5 & 180$^{\circ}$ rotation & 2 & Fig. \ref{fig5}d \\
6 & one diagonal reflection & 2 & Fig. \ref{fig12}d \\
7 & one orthogonal reflection & 2 & Fig. \ref{fig5}b \\
\end{tabular}
\caption{The seven possible symmetries for an English or French board position.}
\label{tab1}
\end{center} 
\end{table}

The English central game is interesting because the board begins and ends
at positions with square symmetry.
John Beasley proved \cite{Beasley} that no matter how the central game is solved, 
the board cannot pass through an intermediate position with square symmetry (type 1)
or 90 degree rotational symmetry (type 2).
We note that this {\it does not mean} that a square symmetric board position
cannot be reached starting with the centre vacant,
only that if a square symmetric position is reached,
it is not solvable.
In what follows we will determine what types of symmetric board positions can appear
during a solution to the central game.

\section{Position class theory}
\label{sec:positionclass}

\begin{figure}[htb]
\centering
\epsfig{file=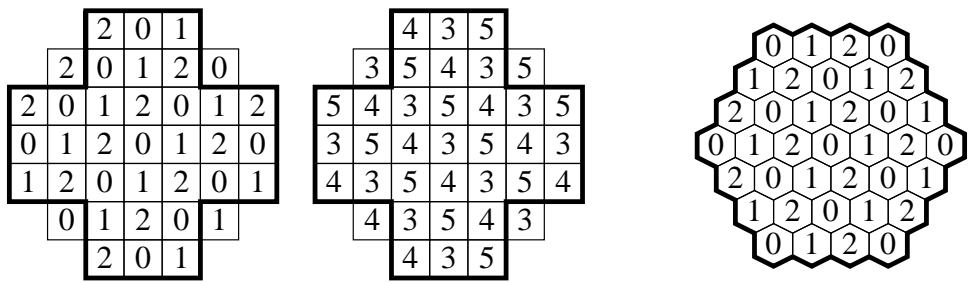}
\caption{Diagonal labeling of holes for the English and French boards (left), and the hexagon board (right).}
\label{fig4}
\end{figure}

Given a board position, if we determine its position class we will know which finishing holes are possible.
We begin with the English and French boards by labeling the holes diagonally with the numbers 0-2,
and again with 3-5, as shown in Figure~\ref{fig4}.

Let $N_i$ be the number of pegs in the holes labeled $i$.
We now observe what happens to $N_0$, $N_1$, $N_2$ after a peg solitaire jump is executed.
One of the three increases by $1$, while the other two decrease by $1$.
Therefore, if we add any two of $N_0$, $N_1$, $N_2$, {\it the parity of the sum can never change as the game is played}.
For example, $(N_1+N_2)\mod 2$ is an invariant of the game, its value can only be $0$ or $1$.
The same holds for $N_3$, $N_4$, $N_5$, so the binary 6-tuple
\begin{equation}
\vec{N} = (N_1+N_2, N_0+N_2, N_0+N_1, N_4+N_5, N_3+N_5, N_3+N_4)
\end{equation}
is an invariant of the game (here each component of $\vec{N}$ is taken modulo 2).

The sixteen values of $\vec{N}$ separate all board positions into sixteen equivalence classes \cite{Beasley},
which we call {\bf position classes}.
A peg solitaire game is played entirely in one position class,
so it is interesting to figure out which position classes have representatives with one peg.
The position class of one peg in the centre is an important one and we call it
``position class A", it corresponds to $\vec{N}=(0,1,1,0,1,1)$.
The reader can check that all three board positions in Figure~\ref{fig1} are in position
class A\footnote{For (a), one should find $N_0=N_3=10$, $N_1=N_2=N_4=N_5=11$,
for (b), $N_0=2$, $N_1=N_2=5$, $N_3=6$, $N_4=N_5=3$, and (c), $N_0=1$, $N_1=N_2=6$.}.

\begin{figure}[htb]
\centering
\epsfig{file=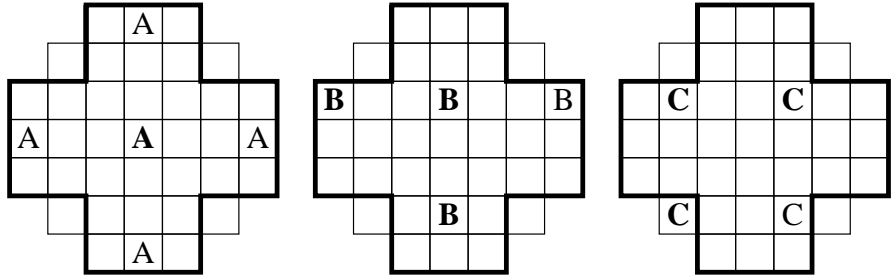}
\caption{The three possible patterns for finishing holes on the English or French boards.
Note that $\vec{N}=(0,1,1,0,1,1)$, $(1,1,0,1,1,0)$ and $(0,1,1,1,0,1)$, respectively.}
\label{fig2}
\end{figure}

Figure~\ref{fig2} shows the possible finishing holes for any board in position class A,
as well as the other one-peg position classes B and C.
Figure~\ref{fig2} shows that there are essentially only three possible patterns for the finishing peg,
up to rotations and reflections.
It should be noted that the holes labeled `A' in Figure~\ref{fig2}a are only a {\it necessary}
condition for a one-peg finish.
If we are in position class A, the only possible finishing holes are those marked by A's.
However, it may not be possible to finish with one peg at all, or only at some A's.

The symmetry of the three patterns in Figure~\ref{fig2} turns out to be very important for what follows.
We note that the pattern of A's has square symmetry (type 1),
while the B's and C's have one reflection symmetry
(types 7 and 6, respectively).

One position class to be avoided is the position class of the empty board, $\vec{N}=(0,0,0,0,0,0)$.
The reason to avoid it is that this position class has no representatives with one peg,
so \textit{any board in the position class of the empty board is unsolvable}.
For example, if we take either of the 37-hole boards and fill the board, then remove the central peg,
we are in the position class of the empty board and therefore in an unsolvable board position.

The same idea can be applied to boards on a triangular lattice.
On the 37-hole hexagon board we use the hole labeling of Figure~\ref{fig4}c,
the vector $\vec{N}$ has only three components,
and there are only four position classes (for details see \cite{BellSol}).
Three of the four position classes have representatives with one peg,
the exception being the position class of the empty board.

\begin{figure}[htb]
\centering
\epsfig{file=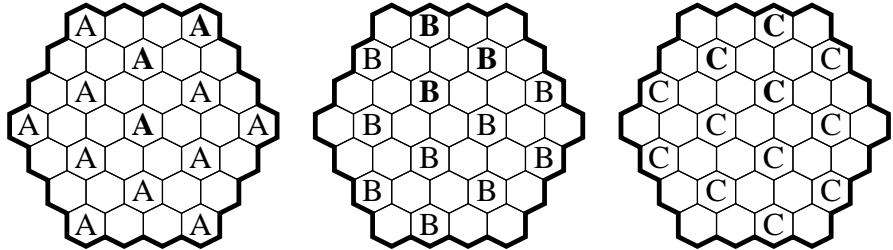}
\caption{The three possible patterns for finishing holes on the hexagon board.}
\label{fig3}
\end{figure}

The three possible patterns for finishing holes are given in Figure~\ref{fig3}.
As before, the position class of one peg in the centre is ``position class A".
We note that patterns B and C are related by a reflection about the $y$-axis,
so there are essentially two patterns for the finishing holes.

\section{Symmetric board positions}
\label{sec:symmetry}

We denote a board position by a lower case letter,
while sets of board positions will be denoted by upper case letters.
If $b$ is a board position then we denote the complement of $b$ by $\overline{b}$,
and if $B$ is a set of board positions,
$\overline{B}$ is the set of complemented board positions
($b\in \overline{B}$ if and only if $\overline{b}\in B$).

\subsection{The English and French boards}
\label{sec:symEngFr}

Suppose we begin from an English or French board position $b$ which is
square symmetric (type 1) and solvable.
Where can the final peg be?
A powerful observation is that
{\it the set of finishing holes must have the same symmetry type as $b$ itself}.
So on one hand we know the set of finishing holes must be square symmetric,
and we also know it can only be a subset of one of the three patterns in Figure~\ref{fig2}.
But only finishing pattern A is square symmetric, so we must be in position class A,
and we can only finish in the holes marked `A'.
Not only that, if the finishing hole is one of those along the edge of the board,
then by reversing the direction of the last jump, we can always finish in the centre.
This same argument works for any symmetry of type 1-5, so we have proved:

\begin{theorem}
On the English and French boards,
if a board position is solvable and has symmetry type 1-5, 
then it lies in position class A and is solvable to the centre.
\label{th1}
\end{theorem}

It is critical in Theorem~\ref{th1} that the board be solvable.
If a board position is not solvable, then the set of finishing holes is empty,
which is trivially a square symmetric pattern.
An unsolvable board may be in the position class of the empty board.
In fact, if we take any board position in position class A,
and remove or add the centre peg (depending on whether it is present or not),
we are in the position class of the empty board,
and therefore not solvable.

What happens if the board position has only a single reflection symmetry (type 6 or 7)?
In that case it may be in position class A, B or C, and it may finish somewhere other than the centre.

Theorem~\ref{th1} tells us that all solvable symmetric positions (types 1-5) are solvable to one peg in the centre.
This suggests that we can find them all by playing backward from one peg in the centre.
Playing peg solitaire backward is equivalent to playing forward from the
complement board position \cite{WinningWays, GPJ04},
so here is a simple algorithm for calculating them:
Let $b_1$ be the board position with a full board with the centre peg missing
(Figure~\ref{fig1}a for the English board),
and define $B_1=\{b_1\}$.
Now define $B_{n+1}$ as the set of all board positions which can be reached from any board position in
$B_n$ by executing any single jump.

Note that by design, $\overline{B}_n$ is the set of $n$-peg board positions solvable to the centre.
We now search through $\overline{B}_n$ for board positions with various symmetries.
The set of all $n$-peg solvable positions with type $j$ symmetry
can be found by searching $\overline{B}_n$ for board positions with type $j$ symmetry.

For details on how these calculations are done, see \cite{BellNotes}.
We do not store duplicate copies of board positions which are rotations or reflections of one another,
each symmetric board position has a single entry, determined by the {\tt mincode()}
(the minimum value of the board code over all symmetry transformations).
A board position is conveniently (but not efficiently) stored in a single, 64-bit integer.

\begin{figure}[htb]
\centering
\epsfig{file=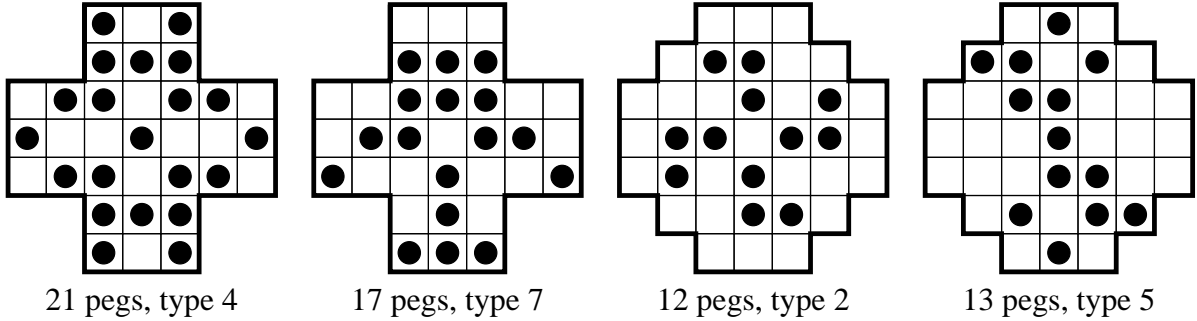}
\caption{Sample solvable boards with an assortment of symmetry types.}
\label{fig5}
\end{figure}

\begin{table}[htb]
\begin{center} 
\begin{tabular}{ c  l  c  c  r  r  r }
& & & Position & English & French & Square \\
Type & Symmetry description & Order & class & 33-hole & 37-hole & 36-hole\\
\hline
1 & square symmetry & 8 & A & 13 & 17 & 21 \\
2 & 90$^{\circ}$ rotation & 4 & A & 25 & 27 & 79 \\
3 & both diagonal reflections & 4 & A & 22 & 126 & 238 \\
4 & both orthogonal reflections & 4 & A & 220 & 258 & 76 \\
5 & 180$^{\circ}$ rotation & 2 & A & 2,238 & 7,051 & 9,148 \\
6 & one diagonal reflection & 2 & A & 5,139 & 40,722 & 64,135 \\
6 & one diagonal reflection & 2 & C & 15,187 & n/c & n/c \\
7 & one orthogonal reflection & 2 & A & 34,501 & 113,375 & 20,961 \\
7 & one orthogonal reflection & 2 & B & 92,732 & n/c & n/c \\
\hline
 & Total & & & 150,077 & 161,576 & 94,658 \\
\end{tabular}
\caption{A count of solvable board positions for the various symmetry types.
``n/c" means not calculated.} 
\label{tab2}
\end{center} 
\end{table}

Table~\ref{tab2} shows the results of such a computation, and four sample positions are shown in Figure~\ref{fig5}.
For the English board, the totals for type 1 and 2 symmetries have been calculated by Beasley \cite{Beasley},
and our results agree with his.
We note that any board position appearing in the English list is solvable on the
French board as well.
We have removed these duplicates,
so the 17 type 1 positions on the French board do not include the 13 English board positions.

The largest $B_i$ for the English board has size $|B_{18}| = 3,626,632$ and for the French board,
$|B_{20}| = 53,371,113$.
This is small enough that a binary search tree of these sets fits into
memory\footnote{Run on a PC with a clock speed of 2.4 GHz and 8 GB of RAM.}.

By Theorem~\ref{th1}, for symmetry type 1-5 we need only start in the centre and all board positions are in position class A.
For symmetry type 6, we need a separate run for position class C, and for type 7, position class B.
Note that for these extra runs, the starting set $B_1$ contains more than one board position.
For example, for position class B we begin with a full board with one peg missing at each B in Figure~\ref{fig3}b,
although due to symmetry it suffices to use only those in bold.
The reason for this is that we need to capture all possible finishing holes.
These runs are time consuming for the French board, and we have not completed them.

Figure~\ref{fig5}c is an interesting case, because this board position
fits on the English board but is not solvable there.
The four added holes are necessary in order to solve it.
You can try all type 1-4 puzzles on my Javascript web program for the English board \cite{BellSymEng}
and the French Board \cite{BellSymFr} (solutions can also be displayed).
You can make these puzzles more challenging by trying to solve them in the minimum number of moves
(where a {\bf move} is one or more jumps by the same peg).

Having computed these symmetrical board positions, we are in a position to demonstrate:

\begin{theorem}
A solution to the central game on the English board cannot pass through an intermediate
position with symmetry type 1-5.
\label{th2}
\end{theorem}

Proof: A board position $b$ can appear during a solution to the central game if and only if
$b$ is solvable to the centre
and $\overline{b}$ is solvable to the centre \cite{GPJ04}.
If $S_j$ is the set of all solvable board positions with type $j$ symmetry,
then a board position $b$ with type $j$ symmetry
can appear during the central game if and only if $b\in S_j$ and $\overline{b}\in S_j$.
We can easily check each element of $S_j$, and we find no matches among types 1-5
(except, of course, for the initial and final board states).

Beasley \cite{Beasley} proves Theorem~\ref{th2} for symmetry types 1 and 2.
In \cite{Beasley180} he proves Theorem~\ref{th2} for type 5 symmetry (180$^{\circ}$ rotation).

Among the $5,139$ type 6 board positions (in position class A),
we find $198$ which form $99$ complement pairs, all can appear during a solution to the central game.
Similarly, there are $912$ type 7 board positions which form 456 complement pairs.
Martin Gardner gave a solution he calls {\bf Jabberwocky} \cite{Gardner} which passes through
eleven intermediate positions with reflection symmetry about the $y$-axis (type 7).
On my web site \cite{BellWeb} I show a solution to the central game which passes through seven
positions with diagonal reflection symmetry (type 6).

\subsection{The 36-hole square board}
\label{sec:sym36}

\begin{figure}[htb]
\centering
\epsfig{file=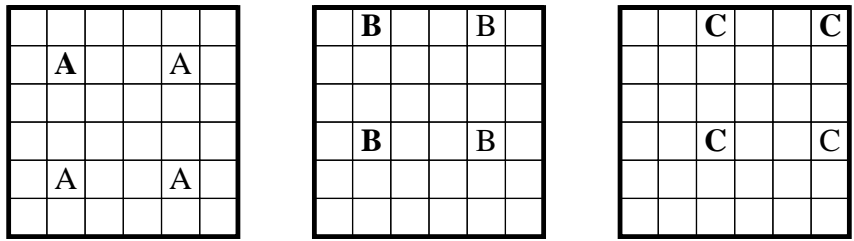}
\caption{The three possible patterns for finishing holes on the square $6\times 6$ board.}
\label{fig11}
\end{figure}

\begin{figure}[htb]
\centering
\epsfig{file=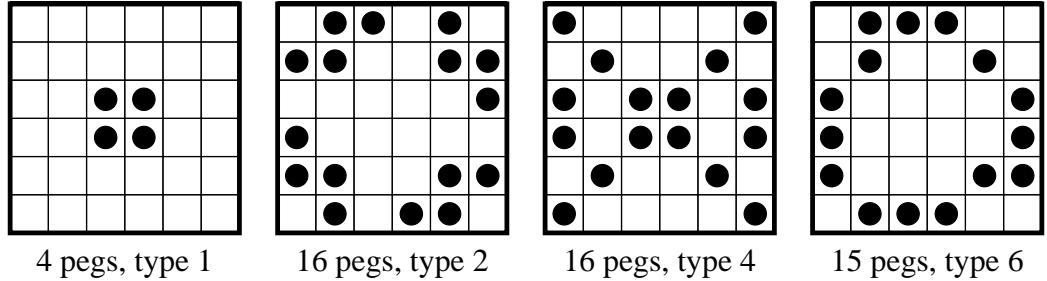}
\caption{Sample solvable $6\times 6$ boards with various symmetry types.}
\label{fig12}
\end{figure}

This board is different because it does not have a central hole.
Nonetheless, it can be analyzed for symmetrical board positions using the same technique.
The smallest solvable board position with square symmetry is
``four pegs in the centre",
shown in Figure~\ref{fig12}a, this board position defines ``position class A".
Figure~\ref{fig11} shows the three possible patterns for finishing holes,
crucially their symmetry types are the same as those in Figure~\ref{fig2}.

Theorem~\ref{th1} is valid on this board
(except for the part about being ``solvable to the centre").
Symmetrical board positions may be calculated by playing backwards,
and Table~\ref{tab2} (right column) includes totals for each symmetry type.
Four examples of symmetrical board positions on the 36-hole square board 
are shown in Figure~\ref{fig12}.
Unfortunately, most symmetrical board positions are easy to solve,
because an obvious sequence of symmetrical jumps reduces them to ``four pegs in the centre".
Figure~\ref{fig12}b-d show three of the harder examples where this is not possible.

\subsection{The 37-hole hexagon board}
\label{sec:symHex}

We now repeat the analysis of symmetrical board positions for the 37-hole hexagon board.
We need to be aware of several important differences.
First, the symmetries are subgroups of $D_{12}$, the dihedral group of order 12
(the symmetries of the regular hexagon).
There are nine possible symmetries, shown in Table~\ref{tab3}.

\begin{table}[htb]
\begin{center} 
\begin{tabular}{ c  l  c  c  r  l }
 & & & Position & & \\
Type & Symmetry description & Order & class & Count & Examples\\
\hline
1 & hexagonal symmetry & 12 & A & 20 & Fig. \ref{fig6}a \\
2 & 60$^{\circ}$ rotation & 6 & A & 14 & Fig. \ref{fig1}c \\
3 & 120$^{\circ}$ rotation \& diagonal refl. ($y$-axis) & 6 & A & 30 & Fig. \ref{fig6}b \\
4 & 120$^{\circ}$ rotation \& orthogonal refl. ($x$-axis) & 6 & A & 87 & Fig. \ref{fig6}c \\
4 & 120$^{\circ}$ rotation \& orthogonal refl. ($x$-axis) & 6 & B & 185 & Fig. \ref{fig6}d \\
5 & both $x$-axis and $y$-axis reflections & 4 & A & 1,438 & Fig. \ref{fig10}c \\
6 & 120$^{\circ}$ rotation & 3 & A & 330 & \\
6 & 120$^{\circ}$ rotation & 3 & B & 754 & Fig. \ref{fig10}b \\
7 & 180$^{\circ}$ rotation & 2 & A & 34,894 & \\
8 & one diagonal reflection ($y$-axis) & 2 & A & 219,295 & Fig. \ref{fig10}d \\
9 & one orthogonal reflection ($x$-axis) & 2 & A & 436,697 & \\
9 & one orthogonal reflection ($x$-axis) & 2 & B & n/c & \\
\end{tabular}
\caption{The nine possible symmetries for a board position on the 37-hole hexagon.
``n/c" means not calculated.} 
\label{tab3}
\end{center} 
\end{table}

We now consider the symmetry of the two possible patterns for finishing pegs on the board,
shown in Figure~\ref{fig3}.  Position class A has hexagonal symmetry, and position class B (and C) have
type 4 symmetry (120$^{\circ}$ rotation plus reflection about the $x$-axis),
as well as the ``sub-symmetries" type 6 and 9.
The same argument as before leads to:

\begin{theorem}
On the 37-hole hexagon board,
if a board position is solvable and has symmetry type 1-3, 5, 7 or 8, 
then it lies in position class A.
\label{th3}
\end{theorem}

Another difference is that when the board is solvable and in position class A,
it may not be solvable to the centre (as in Figure~\ref{fig1}c).
Therefore, when we initialize the set $B_1$,
we need to start with three board positions with one peg missing at each of the
bold A's in Figure~\ref{fig3}a.
The calculation of all $B_n$ for position class A is time consuming,
taking about a week of CPU time and 20 GB of disk space.
The largest set $B_n$ in this case is $|B_{19}| = 364,696,466$,
and the binary search tree containing it is too large to fit into memory on my machine.
The calculation therefore has to be split into pieces, which increases the computation time.

\begin{figure}[htb]
\centering
\epsfig{file=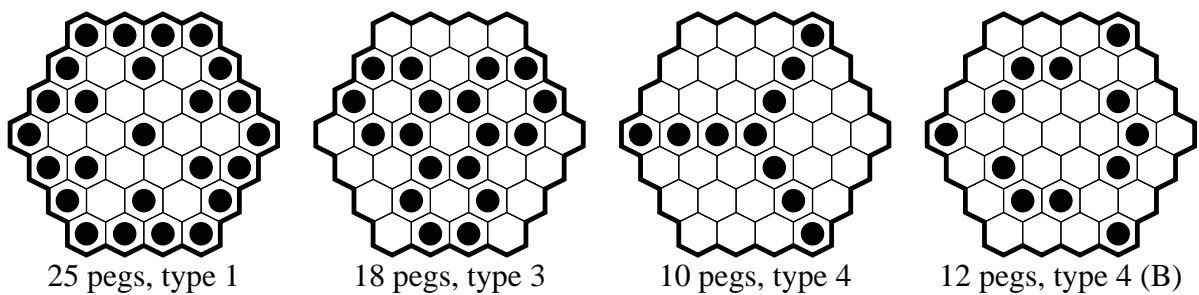}
\caption{Sample solvable hexagon boards with an assortment of symmetry types.}
\label{fig6}
\end{figure}

\begin{figure}[htb]
\centering
\epsfig{file=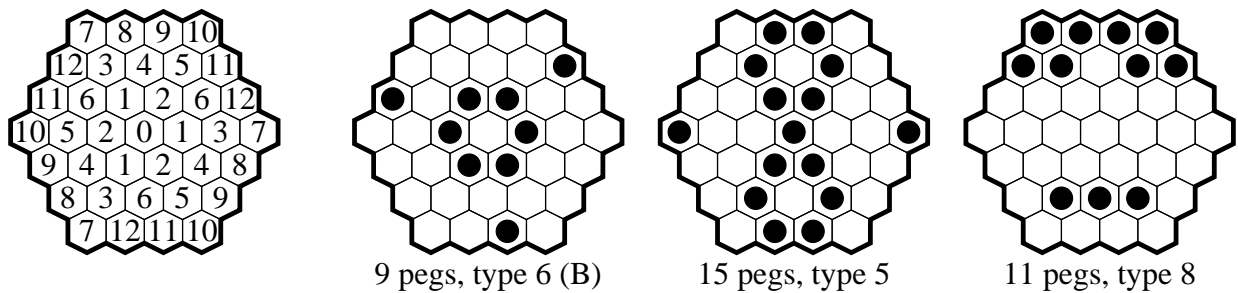}
\caption{A template for 120$^{\circ}$ symmetry, and more sample solvable hexagon boards.}
\label{fig10}
\end{figure}

Figures~\ref{fig6} and \ref{fig10} show seven board positions obtained from these calculations,
with board counts given in Table~\ref{tab3} and Javascript web program for types 1-4 here \cite{BellSymHex}.
Note that Figure~\ref{fig6}d and Figure~\ref{fig10}b show board positions in position class B.
All other board positions shown in this document are in position class A.
The board positions in Figure~\ref{fig10}c and d were selected because they are particularly difficult to solve.

To complete all entries in Table~\ref{tab3} for position class B would require a second run,
even more time consuming than the first.
For position class B and symmetry types 4 and 6,
we used a different technique.
All board positions with 120$^{\circ}$ symmetry can be obtained by mapping every 13-bit binary integer
to the board in Figure~\ref{fig10}a, where a peg is present at location $i$ iff the $i$'th bit is set.
We can then exhaustively attempt to solve each board, one by one,
to derive a complete list of solvable boards by position class and symmetry type (1-4 or 6).
This technique is reasonable when the total number of boards is under a few thousand,
and it gives us a way to double check our results (at least for types 1-4 and 6).

\section{Board positions with a unique winning jump}
\label{sec:unique}

Many of the symmetrical board positions found in the previous section tend to be easy to solve by hand.
The problems shown in Figures~\ref{fig5}-\ref{fig10} are not typical, they are some of the harder problems.
Often, it is possible to make a few symmetrical jumps which reduce the pattern to a smaller symmetrical pattern
which has been solved previously.

At any initial board position, a number of starting jumps are available.
Often, any starting jump can be executed, after which the board remains solvable.
But suppose we search specifically for initial positions where only {\bf one} of the starting jumps
gives a solvable board position.
It is not obvious that such board positions exist,
because during the English central game (for example)
most board positions have many possible jumps that can lead to a solution.

Fortunately, the sets $B_n$ calculated in the last section are exactly what we need to search for these
``unique winning jump" board positions.
Consider a board position $b\in \overline{B}_n$.
We are looking for $b$'s where only a single jump ends at a solvable board position.
We can execute jump $k$ on $b$, producing the board position $b_k$.
$b_k$ is solvable if and only if $b_k\in \overline{B}_{n-1}$.

In order to find a puzzle which is ``difficult",
the number of dead ends should be large,
this suggests we want a large number of starting jumps.
It would seem that the most difficult $n$-peg initial positions are those which
\begin{packed_enumerate}
\item{Have a single winning jump, and}
\item{Have as many starting jumps as possible.}
\end{packed_enumerate}

\begin{figure}[htb]
\centering
\epsfig{file=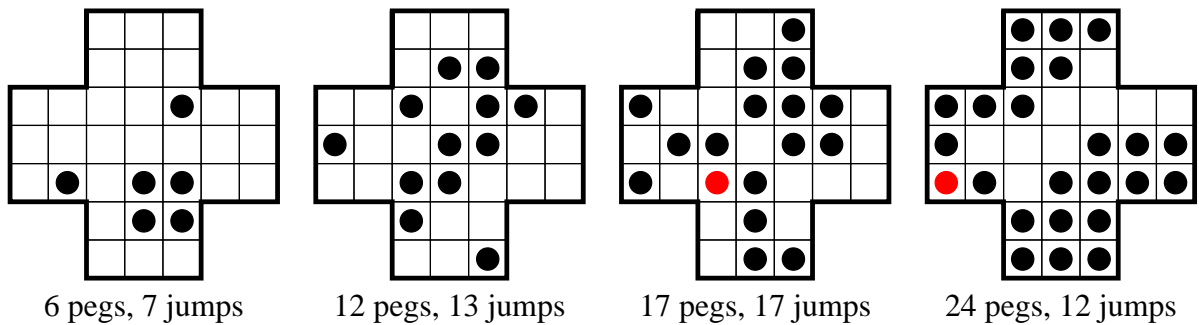}
\caption{English puzzles with a unique winning jump.
Under each diagram is the number of pegs and the total number of starting jumps.
Playable on the web at \cite{BellDifEng}.}
\label{fig7}
\end{figure}
\begin{figure}[htb]
\centering
\epsfig{file=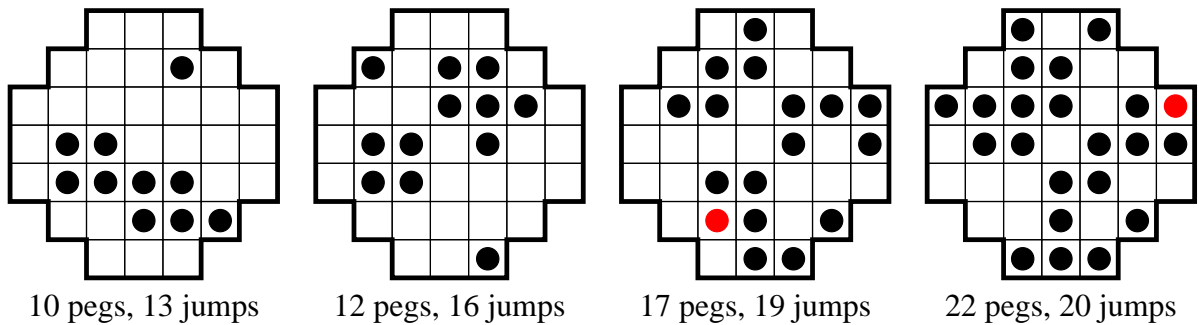}
\caption{French puzzles with a unique winning jump.
Playable on the web at \cite{BellDifFr}.}
\label{fig8}
\end{figure}
\begin{figure}[htb]
\centering
\epsfig{file=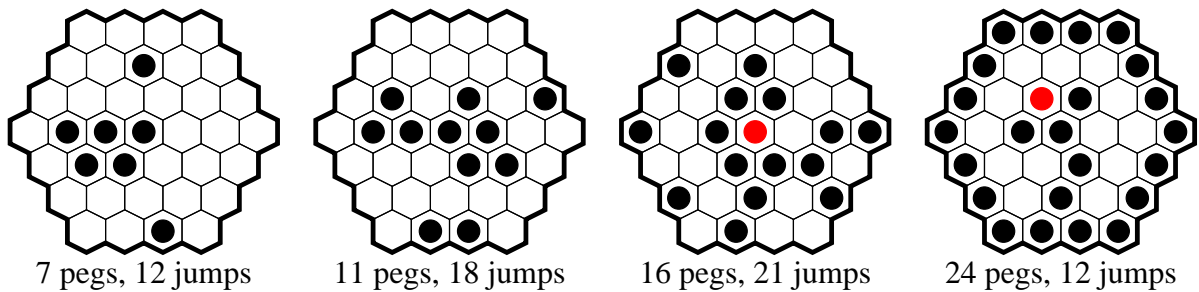}
\caption{Hexagonal puzzles with a unique winning jump.
Playable on the web at \cite{BellDifHex}.}
\label{fig9}
\end{figure}

Table~\ref{tab4} summarizes the results of these calculations,
and Figures~\ref{fig7}-\ref{fig9} show example board positions calculated using this strategy.
All of these puzzles can be played on my Javascript programs \cite{BellDifEng,BellDifFr,BellDifHex}
(the programs can also display solutions).
We note that for a particular board and number of pegs $n$,
there is sometimes a \textit{unique} board position with as many jumps as possible and one winning jump.
When an entry in Table~\ref{tab4} is blank,
this indicates there are no $n$-peg board positions with a unique winning jump.

\begin{table}[!t]
\begin{center} 
\begin{tabular}{ c | c  c | c  c | c  c}
 & \multicolumn{2}{|c|}{English board} & \multicolumn{2}{|c|}{French board} & \multicolumn{2}{|c}{hexagon board}\\
$n$ (pegs) & max jumps & count & max jumps & count & max jumps & count\\
\hline
4 & 4 & 2 & 4 & 3 & 6 & 6\\ 
5 & 7 & 1 & 7 & 1 & 10 & 2\\
6 & 7 & 2$\dagger$ & 8 & 1 & 10 & 26\\
7 & 9 & 2 & 10 & 1 & 12 & 10$\dagger$\\
8 & 9 & 4 & 10 & 1 & 14 & 6\\
9 & 11 & 1 & 12 & 1 & 16 & 2\\
10 & 12 & 1 & 13 & 1$\dagger$ & 18 & 1\\
11 & 12 & 2 & 14 & 2 & 18 & 3$\dagger$\\
12 & 13 & 4$\dagger$ & 16 & 1$\dagger$ & 20 & 1\\
13 & 14 & 1 & 17 & 1 & 20 & 1\\
14 & 14 & 4 & 16 & 4 & 22 & 1\\
15 & 15 & 1 & 17 & 2 & 21 & 2\\
16 & 15 & 1 & 19 & 1 & 21 & 2$\dagger$\\
17 & 17 & 1$\dagger$ & 19 & 2$\dagger$ & 22 & 1\\
18 & 15 & 2 & 18 & 4 & 22 & 1\\
19 & 15 & 2 & 21 & 1 & 20 & 1\\
20 & 14 & 3 & 19 & 1 & 18 & 3\\
21 & 13 & 3 & 18 & 5 & 17 & 1\\
22 & 14 & 1 & 20 & 1$\dagger$ & 16 & 2\\
23 & 11 & 6 & 18 & 1 & 13 & 1\\
24 & 12 & 1$\dagger$ & 18 & 2 & 12 & 1$\dagger$\\
25 & 10 & 1 & 17 & 1 & & \\
26 & 7 & 3 & 17 & 1 & & \\
27 & 6 & 2 & 16 & 1 & & \\
28 & & & 13 & 1 & & \\
29 & & & 11 & 1 & & \\
30 & & & 10 & 3 & & \\
31 & & & 6 & 1 & & \\
32 & & & 8 & 1 & & \\
\end{tabular}
\caption{A summary of board positions with a unique winning jump by pegs and maximum starting jumps,
for each of the three board types.
$\dagger$: case appears in Figures~\ref{fig7}-\ref{fig9}.}.
\label{tab4}
\end{center} 
\end{table}

These puzzles tend to be challenging to solve by hand, particularly as they become larger.
As an aid to the solver, we have identified the first peg to jump in red for the larger board positions.

If some jump is a winner on a symmetrical board position, then the symmetrical
partners of this jump are also winners.
Thus, board positions with a unique winning jump tend to have no symmetry.
The only exception would be a board position with a single reflection symmetry,
it could have a single winning jump along the axis of reflection.
We have found a few examples like this, but none with the maximum number of jumps.

These puzzles have an entirely different character from the symmetrical puzzles in the last section.
The fact that there is a unique winning jump can be used to help solve these puzzles.
After the first jump is executed, any jump which could have been executed first must still be a dead end.
The second jump can only be a jump which was opened up by the first jump,
and so on.
Sometimes, if you can identify the first jump, the rest of the solution follows more easily.

While the first jump of a solution is unique,
subsequent jumps can often be executed in either order,
or the final jump can go in either direction, so the solution is not unique.
However, a few of these puzzles do have a unique solution, which is quite rare in peg solitaire.
An example is Figure~\ref{fig9}c---there is only one sequence of jumps which solves this puzzle.

When the number of pegs is relatively small (say, under 13), the board may not limit the jumping possibilities.
We can often translate the pattern of pegs, and this gives a solution which is counted as different.
This effect is responsible for the large counts on the hexagonal board (26 solutions with 6 pegs and 10 jumps),
this is not 26 different solutions but a few solutions translated.
These board positions with less than 13 pegs can be considered as puzzles on an infinite board.
These puzzles retain the property that the number of starting jumps is large,
and there is a unique winning jump.
As the puzzles become larger, on an infinite board the unique winning jump property tends to be lost.

Finally, we note that all these unique winning jump puzzles were calculated using position class A.
Since symmetry plays no role here, there is no reason why we could not use position class B or C.
This would produce a whole new set of problems with a unique winning jump,
and Table~\ref{tab4} would be different for each position class.

\section{Summary}
\label{sec:summary}

We have presented two different strategies for creating peg solitaire puzzles.
The first searches for solvable symmetric positions,
while the second identifies solvable positions with a unique winning jump.
The two strategies don't seem to have anything in common,
but they can both be calculated using the sets of board positions $B_n$
obtained by playing the game backwards.

The central game on the English board is an attractive puzzle because
it begins and ends at positions with square symmetry,
but in between symmetry is lost, and by Theorem~\ref{th2} symmetry is not possible
(except for reflection symmetry).
For a good puzzle,
it is desirable that symmetry is not possible in the middle,
for symmetric intermediate positions often indicate an easy solution where jumps are repeated in a symmetrical fashion.
We have identified all solvable symmetric board positions (of most types), both on a square a triangular grid.
Many of these make nice puzzles to solve by hand.

The ``unique winning jump" puzzles have a completely different feel---they lack
symmetry and are much harder to solve.
The fact that they have a unique winning jump can be exploited by crafty solvers.

Any solvable board position presented above is also solvable when considered on an infinite board.
This means that in some sense these puzzles exist independently of any particular board.
We saw in going from the English to French board that additional puzzles were found that were solvable on the
French board but not on the English board (Figure~\ref{fig5}c).
Similarly, in going from the French board to an infinite board,
we would expect additional problems solvable only on a sufficiently large board.
Searching for all $n$-peg symmetric or unique winning jump puzzles
on an infinite board is an interesting computational challenge.


\begin{thebibliography}{mybib} 

\bibitem{Beasley} J. Beasley, \textit{The Ins and Outs of Peg Solitaire}, Oxford Univ. Press, 1992.

\bibitem{WinningWays} E. Berlekamp, J. Conway and R. Guy, Purging pegs properly, in \textit{Winning Ways for Your Mathematical Plays},
2nd ed., Vol. 4, Chap. 23: 803--841, A K Peters, 2004.

\bibitem{Gardner} M. Gardner, Peg Solitaire, in \textit{Knots and Borromean Rings, Rep-Tiles and Eight Queens},
Cambridge Univ. Press, 2014

\bibitem{GPJ04} G. Bell, Triangular peg solitaire unlimited, {\it Games and Puzzles J.} \#36 (2004),\newline
\href{http://www.gpj.connectfree.co.uk/gpjr.htm}{\tt http://www.gpj.connectfree.co.uk/gpjr.htm}\newline
\href{http://arxiv.org/abs/0711.0486}{\tt http://arxiv.org/abs/0711.0486}

\bibitem{BellSol} G. Bell, Solving triangular peg solitaire, \textit{J. Integer Sequences}, 08.4.8, \textbf{11} (2008),\newline
\href{http://arxiv.org/abs/math/0703865}{\tt http://arxiv.org/abs/math/0703865}

\bibitem{Beasley180} J. Beasley, On 33-hole solitaire positions with rotational symmetry (2012),\newline
\href{http://www.jsbeasley.co.uk/puzzles/solitairerotations.pdf}
{\tt http://www.jsbeasley.co.uk/puzzles/solitairerotations.pdf}

\bibitem{BellNotes} G. Bell, Notes on solving and playing peg solitaire on a computer (2014), \newline
\href{http://arxiv.org/abs/0903.3696}{\tt http://arxiv.org/abs/0903.3696}

\bibitem{BellWeb} G. Bell, Peg Solitaire web site,
\href{http://www.gibell.net/pegsolitaire/}
{\tt http://www.gibell.net/pegsolitaire/}

\bibitem{BellSymEng} G. Bell, Symmetric English positions,\newline
\href{http://www.gibell.net/pegsolitaire/Tools/Symmetric/English.htm}
{\tt http://www.gibell.net/pegsolitaire/Tools/Symmetric/English.htm}

\bibitem{BellSymFr} G. Bell, Symmetric French positions,\newline
\href{http://www.gibell.net/pegsolitaire/Tools/Symmetric/French.htm}
{\tt http://www.gibell.net/pegsolitaire/Tools/Symmetric/French.htm}

\bibitem{BellSym6x6} G. Bell, Symmetric 6x6 positions,\newline
\href{http://www.gibell.net/pegsolitaire/Tools/Symmetric/SixBySix.htm}
{\tt http://www.gibell.net/pegsolitaire/Tools/Symmetric/SixBySix.htm}

\bibitem{BellSymHex} G. Bell, Symmetric Hexagon positions,\newline
\href{http://www.gibell.net/pegsolitaire/Tools/Hex37/Symmetric.htm}
{\tt http://www.gibell.net/pegsolitaire/Tools/Hex37/Symmetric.htm}

\bibitem{BellDifEng} G. Bell, Difficult English positions,\newline
\href{http://www.gibell.net/pegsolitaire/Tools/Difficult/English0.htm}
{\tt http://www.gibell.net/pegsolitaire/Tools/Difficult/English0.htm}

\bibitem{BellDifFr} G. Bell, Difficult French positions,\newline
\href{http://www.gibell.net/pegsolitaire/Tools/Difficult/French0.htm}
{\tt http://www.gibell.net/pegsolitaire/Tools/Difficult/French0.htm}

\bibitem{BellDifHex} G. Bell, Difficult Hexagon positions,\newline
\href{http://www.gibell.net/pegsolitaire/Tools/Hex37/Difficult.htm}
{\tt http://www.gibell.net/pegsolitaire/Tools/Hex37/Difficult.htm}

\end{thebibliography}
\end{document}